\documentclass{amsart}
\usepackage{amsfonts}

\newtheorem{theorem}{Theorem}
\theoremstyle{plain}

\newtheorem{corollary}{Corollary}

\newtheorem{proposition}{Proposition}
\newtheorem{remark}{Remark}

\numberwithin{equation}{section}
\input{tcilatex}

\begin{document}
\title[Reverses of the Continuous Triangle Inequality]{Some Reverses of the
Continuous Triangle Inequality for Bochner Integral of Vector-Valued
Functions in Complex Hilbert Spaces}
\author{Sever S. Dragomir}
\address{School of Computer Science and Mathematics\\
Victoria University of Technology\\
PO Box 14428, MCMC 8001\\
Victoria, Australia.}
\email{sever@csm.vu.edu.au}
\urladdr{http://rgmia.vu.edu.au/SSDragomirWeb.html}
\date{April 16, 2004.}
\subjclass[2000]{Primary 46C05; Secondary 26D15.}
\keywords{Triangle inequality, Diaz-Metcalf inequality, Integral
inequalities, Complex Hilbert spaces.}

\begin{abstract}
Some reverses of the continuous triangle inequality for Bochner integral of
vector-valued functions in complex Hilbert spaces are given. Applications
for complex-valued functions are provided as well.
\end{abstract}

\maketitle

\section{Introduction}

Let $f:\left[ a,b\right] \rightarrow \mathbb{K}$, $\mathbb{K}=\mathbb{C}$ or 
$\mathbb{R}$ be a Lebesgue integrable function. The following inequality is
the continuous version of the triangle inequality%
\begin{equation}
\left\vert \int_{a}^{b}f\left( x\right) dx\right\vert \leq
\int_{a}^{b}\left\vert f\left( x\right) \right\vert dx,  \label{1.1}
\end{equation}%
and plays a fundamental role in Mathematical Analysis and its applications.

It appears, see \cite[p. 492]{MPF}, that the first reverse inequality for (%
\ref{1.1}) was obtained by J. Karamata in his book from 1949, \cite{K}:%
\begin{equation}
\cos \theta \int_{a}^{b}\left\vert f\left( x\right) \right\vert dx\leq
\left\vert \int_{a}^{b}f\left( x\right) dx\right\vert  \label{1.2}
\end{equation}%
provided%
\begin{equation*}
-\theta \leq \arg \left[ f\left( x\right) \right] \leq \theta ,\ \ x\in %
\left[ a,b\right]
\end{equation*}%
for given $\theta \in \left( 0,\frac{\pi }{2}\right) .$

In \cite{SSD1}, the author has extended the above result for Bochner
integrals of vector-valued functions in real or complex Hilbert spaces.

If $\left( H;\left\langle \cdot ,\cdot \right\rangle \right) $ is a Hilbert
space over $\mathbb{K}$ $\left( \mathbb{K}=\mathbb{C},\mathbb{R}\right) $
and $f\in L\left( \left[ a,b\right] ;H\right) $, this means that $f:\left[
a,b\right] \rightarrow H$ is Bochner measurable on $\left[ a,b\right] $ and $%
\int_{a}^{b}\left\Vert f\left( t\right) \right\Vert dt$ is finite, and there
exists a constant $K\geq 1$ and a vector $e\in H,$ $\left\Vert e\right\Vert
=1$ such that%
\begin{equation}
\left\Vert f\left( t\right) \right\Vert \leq K\func{Re}\left\langle f\left(
t\right) ,e\right\rangle \text{ \ for a.e. \ }t\in \left[ a,b\right] ,
\label{1.3}
\end{equation}%
then we have the inequality%
\begin{equation}
\int_{a}^{b}\left\Vert f\left( t\right) \right\Vert dt\leq K\left\Vert
\int_{a}^{b}f\left( t\right) dt\right\Vert .  \label{1.4}
\end{equation}%
This provides a reverse inequality for the well know result for Bochner
integrals and vector-valued functions:%
\begin{equation}
\left\Vert \int_{a}^{b}f\left( t\right) dt\right\Vert \leq
\int_{a}^{b}\left\Vert f\left( t\right) \right\Vert dt  \label{1.5}
\end{equation}%
for any $f\in L\left( \left[ a,b\right] ;H\right) .$

Note that the case of equality holds in (\ref{1.4}) (see \cite{SSD1}) if and
only if 
\begin{equation}
\int_{a}^{b}f\left( t\right) dt=\frac{1}{K}\left( \int_{a}^{b}\left\Vert
f\left( t\right) \right\Vert dt\right) e.  \label{1.6}
\end{equation}

For some particular cases of interest, see \cite{SSD1}.

The main aim of the present paper is to point out some newer inequalities
for complex Hilbert spaces under various conditions for both $\func{Re}%
\left\langle f\left( t\right) ,e\right\rangle $ and $\func{Im}\left\langle
f\left( t\right) ,e\right\rangle $ $\left( e\in H,\ \left\Vert e\right\Vert
=1\right) $ and in this way improve some earlier results from \cite{SSD1}
that have been stated for real or complex Hilbert spaces. Applications for
complex-valued functions are also provided.

\section{The Case of a Unit Vector}

The following result holds.

\begin{theorem}
\label{t2.1}Let $\left( H;\left\langle \cdot ,\cdot \right\rangle \right) $
be a complex Hilbert space. If $f\in L\left( \left[ a,b\right] ;H\right) $
is such that there exists $k_{1},k_{2}\geq 0$ with%
\begin{equation}
k_{1}\left\Vert f\left( t\right) \right\Vert \leq \func{Re}\left\langle
f\left( t\right) ,e\right\rangle ,\ \ k_{2}\left\Vert f\left( t\right)
\right\Vert \leq \func{Im}\left\langle f\left( t\right) ,e\right\rangle
\label{2.1}
\end{equation}%
for a.e. $t\in \left[ a,b\right] ,$ where $e\in H,$ $\left\Vert e\right\Vert
=1,$ is given, then%
\begin{equation}
\sqrt{k_{1}^{2}+k_{2}^{2}}\int_{a}^{b}\left\Vert f\left( t\right)
\right\Vert dt\leq \left\Vert \int_{a}^{b}f\left( t\right) dt\right\Vert .
\label{2.2}
\end{equation}%
The case of equality holds in (\ref{2.2}) if and only if%
\begin{equation}
\int_{a}^{b}f\left( t\right) dt=\left( k_{1}+ik_{2}\right) \left(
\int_{a}^{b}\left\Vert f\left( t\right) \right\Vert dt\right) e.  \label{2.3}
\end{equation}
\end{theorem}

\begin{proof}
Using the Schwarz inequality $\left\Vert u\right\Vert \left\Vert
v\right\Vert \geq \left\vert \left\langle u,v\right\rangle \right\vert ,$ $%
u,v\in H;$ in the complex Hilbert space $\left( H;\left\langle \cdot ,\cdot
\right\rangle \right) ,$ we have%
\begin{align}
\left\Vert \int_{a}^{b}f\left( t\right) dt\right\Vert ^{2}& =\left\Vert
\int_{a}^{b}f\left( t\right) dt\right\Vert ^{2}\left\Vert e\right\Vert ^{2}
\label{2.4} \\
& \geq \left\vert \left\langle \int_{a}^{b}f\left( t\right)
dt,e\right\rangle \right\vert ^{2}=\left\vert \int_{a}^{b}\left\langle
f\left( t\right) ,e\right\rangle dt\right\vert ^{2}  \notag \\
& =\left\vert \int_{a}^{b}\func{Re}\left\langle f\left( t\right)
,e\right\rangle dt+i\left( \int_{a}^{b}\func{Im}\left\langle f\left(
t\right) ,e\right\rangle dt\right) \right\vert ^{2}  \notag \\
& =\left( \int_{a}^{b}\func{Re}\left\langle f\left( t\right) ,e\right\rangle
dt\right) ^{2}+\left( \int_{a}^{b}\func{Im}\left\langle f\left( t\right)
,e\right\rangle dt\right) ^{2}.  \notag
\end{align}%
Now, on integrating (\ref{2.1}), we deduce%
\begin{equation}
k_{1}\int_{a}^{b}\left\Vert f\left( t\right) \right\Vert dt\leq \int_{a}^{b}%
\func{Re}\left\langle f\left( t\right) ,e\right\rangle dt,\quad
k_{2}\int_{a}^{b}\left\Vert f\left( t\right) \right\Vert dt\leq \int_{a}^{b}%
\func{Im}\left\langle f\left( t\right) ,e\right\rangle dt  \label{2.5}
\end{equation}%
implying%
\begin{equation}
\left( \int_{a}^{b}\func{Re}\left\langle f\left( t\right) ,e\right\rangle
dt\right) ^{2}\geq k_{1}^{2}\left( \int_{a}^{b}\left\Vert f\left( t\right)
\right\Vert dt\right) ^{2}  \label{2.6}
\end{equation}%
and%
\begin{equation}
\left( \int_{a}^{b}\func{Im}\left\langle f\left( t\right) ,e\right\rangle
dt\right) ^{2}\geq k_{2}^{2}\left( \int_{a}^{b}\left\Vert f\left( t\right)
\right\Vert dt\right) ^{2}.  \label{2.7}
\end{equation}%
If we add (\ref{2.6}) and (\ref{2.7}) and use (\ref{2.4}), we deduce the
desired inequality (\ref{2.2}).

Further, if (\ref{2.3}) holds, then obviously%
\begin{align*}
\left\Vert \int_{a}^{b}f\left( t\right) dt\right\Vert & =\left\vert
k_{1}+ik_{2}\right\vert \left( \int_{a}^{b}\left\Vert f\left( t\right)
\right\Vert dt\right) \left\Vert e\right\Vert \\
& =\sqrt{k_{1}^{2}+k_{2}^{2}}\int_{a}^{b}\left\Vert f\left( t\right)
\right\Vert dt,
\end{align*}%
and the equality case holds in (\ref{2.2}).

Before we prove the reverse implication, let us observe that, for $x\in H$
and $e\in H,$ $\left\Vert e\right\Vert =1,$ the following identity is valid%
\begin{equation*}
\left\Vert x-\left\langle x,e\right\rangle e\right\Vert ^{2}=\left\Vert
x\right\Vert ^{2}-\left\vert \left\langle x,e\right\rangle \right\vert ^{2},
\end{equation*}%
therefore $\left\Vert x\right\Vert =\left\vert \left\langle x,e\right\rangle
\right\vert $ if and only if $x=\left\langle x,e\right\rangle e.$

If we assume that equality holds in (\ref{1.2}), then the case of equality
must hold in all the inequalities required in the argument used to prove the
inequality (\ref{2.2}). Therefore, we must have%
\begin{equation}
\left\Vert \int_{a}^{b}f\left( t\right) dt\right\Vert =\left\vert
\left\langle \int_{a}^{b}f\left( t\right) dt,e\right\rangle \right\vert
\label{2.7a}
\end{equation}%
and 
\begin{equation}
k_{1}\left\Vert f\left( t\right) \right\Vert =\func{Re}\left\langle f\left(
t\right) ,e\right\rangle ,\ \ \ k_{2}\left\Vert f\left( t\right) \right\Vert
=\func{Im}\left\langle f\left( t\right) ,e\right\rangle  \label{2.8}
\end{equation}%
for a.e. $t\in \left[ a,b\right] .$

From (\ref{2.7a}) we deduce%
\begin{equation}
\int_{a}^{b}f\left( t\right) dt=\left\langle \int_{a}^{b}f\left( t\right)
dt,e\right\rangle e,  \label{2.9}
\end{equation}%
and from (\ref{2.8}), by multiplying the second equality with $i,$the
imaginary unit, and integrating both equations on $\left[ a,b\right] ,$ we
deduce%
\begin{equation}
\left( k_{1}+ik_{2}\right) \int_{a}^{b}\left\Vert f\left( t\right)
\right\Vert dt=\left\langle \int_{a}^{b}f\left( t\right) dt,e\right\rangle .
\label{2.10}
\end{equation}%
Finally, by (\ref{2.9}) and (\ref{2.10}), we deduce the desired equality (%
\ref{2.3}).
\end{proof}

The following corollary is of interest.

\begin{corollary}
\label{c2.2}Let $e$ be a unit vector in the complex Hilbert space $\left(
H;\left\langle \cdot ,\cdot \right\rangle \right) $ and $\eta _{1},\eta
_{2}\in \left( 0,1\right) .$ If $f\in L\left( \left[ a,b\right] ;H\right) $
is such that%
\begin{equation}
\left\Vert f\left( t\right) -e\right\Vert \leq \eta _{1},\quad \left\Vert
f\left( t\right) -ie\right\Vert \leq \eta _{2}\text{ \ for a.e. \ }t\in %
\left[ a,b\right] ,  \label{2.11}
\end{equation}%
then we have the inequality%
\begin{equation}
\sqrt{2-\eta _{1}^{2}-\eta _{2}^{2}}\int_{a}^{b}\left\Vert f\left( t\right)
\right\Vert dt\leq \left\Vert \int_{a}^{b}f\left( t\right) dt\right\Vert .
\label{2.12}
\end{equation}%
The case of equality holds in (\ref{2.12}) if and only if%
\begin{equation}
\int_{a}^{b}f\left( t\right) dt=\left( \sqrt{1-\eta _{1}^{2}}+i\sqrt{1-\eta
_{2}^{2}}\right) \left( \int_{a}^{b}\left\Vert f\left( t\right) \right\Vert
dt\right) e.  \label{2.13}
\end{equation}
\end{corollary}

\begin{proof}
From the first inequality in (\ref{2.11}) we deduce, by taking the square,
that%
\begin{equation*}
\left\Vert f\left( t\right) \right\Vert ^{2}+1-\eta _{1}^{2}\leq 2\func{Re}%
\left\langle f\left( t\right) ,e\right\rangle ,
\end{equation*}%
implying%
\begin{equation}
\frac{\left\Vert f\left( t\right) \right\Vert ^{2}}{\sqrt{1-\eta _{1}^{2}}}+%
\sqrt{1-\eta _{1}^{2}}\leq \frac{2\func{Re}\left\langle f\left( t\right)
,e\right\rangle }{\sqrt{1-\eta _{1}^{2}}}  \label{2.14}
\end{equation}%
for a.e. $t\in \left[ a,b\right] .$

Since, obviously%
\begin{equation}
2\left\Vert f\left( t\right) \right\Vert \leq \frac{\left\Vert f\left(
t\right) \right\Vert ^{2}}{\sqrt{1-\eta _{1}^{2}}}+\sqrt{1-\eta _{1}^{2}},
\label{2.15}
\end{equation}%
hence, by (\ref{2.14}) and (\ref{2.15}) we get%
\begin{equation}
0\leq \sqrt{1-\eta _{1}^{2}}\left\Vert f\left( t\right) \right\Vert \leq 
\func{Re}\left\langle f\left( t\right) ,e\right\rangle  \label{2.16}
\end{equation}%
for a.e. $t\in \left[ a,b\right] .$

From the second inequality in (\ref{2.11}) we deduce%
\begin{equation*}
0\leq \sqrt{1-\eta _{2}^{2}}\left\Vert f\left( t\right) \right\Vert \leq 
\func{Re}\left\langle f\left( t\right) ,ie\right\rangle
\end{equation*}%
for a.e. $t\in \left[ a,b\right] .$ Since%
\begin{equation*}
\func{Re}\left\langle f\left( t\right) ,ie\right\rangle =\func{Im}%
\left\langle f\left( t\right) ,e\right\rangle
\end{equation*}%
hence%
\begin{equation}
0\leq \sqrt{1-\eta _{2}^{2}}\left\Vert f\left( t\right) \right\Vert \leq 
\func{Im}\left\langle f\left( t\right) ,e\right\rangle  \label{2.17}
\end{equation}%
for a.e. $t\in \left[ a,b\right] .$

Now, observe from (\ref{2.16}) and (\ref{2.17}), that the condition (\ref%
{2.1}) of Theorem \ref{t2.1} is satisfied for $k_{1}=\sqrt{1-\eta _{1}^{2}},$
$k_{2}=\sqrt{1-\eta _{2}^{2}}\in \left( 0,1\right) ,$ and thus the corollary
is proved.
\end{proof}

The following corollary may be stated as well.

\begin{corollary}
\label{c2.3}Let $e$ be a unit vector in the complex Hilbert space $\left(
H;\left\langle \cdot ,\cdot \right\rangle \right) $ and $M_{1}\geq m_{1}>0,$ 
$M_{2}\geq m_{2}>0.$ If $f\in L\left( \left[ a,b\right] ;H\right) $ is such
that either%
\begin{equation}
\func{Re}\left\langle M_{1}e-f\left( t\right) ,f\left( t\right)
-m_{1}e\right\rangle \geq 0,\ \ \ \text{ \ }\func{Re}\left\langle
M_{2}ie-f\left( t\right) ,f\left( t\right) -m_{2}ie\right\rangle \geq 0
\label{2.18}
\end{equation}%
or, equivalently,%
\begin{align}
\left\Vert f\left( t\right) -\frac{M_{1}+m_{1}}{2}e\right\Vert & \leq \frac{1%
}{2}\left( M_{1}-m_{1}\right) ,  \label{2.19} \\
\left\Vert f\left( t\right) -\frac{M_{2}+m_{2}}{2}ie\right\Vert & \leq \frac{%
1}{2}\left( M_{2}-m_{2}\right) ,  \notag
\end{align}%
for each a.e. $t\in \left[ a,b\right] ,$ then we have the inequality%
\begin{equation}
2\left[ \frac{m_{1}M_{1}}{\left( M_{1}+m_{1}\right) ^{2}}+\frac{m_{2}M_{2}}{%
\left( M_{2}+m_{2}\right) ^{2}}\right] ^{1/2}\int_{a}^{b}\left\Vert f\left(
t\right) \right\Vert dt\leq \left\Vert \int_{a}^{b}f\left( t\right)
dt\right\Vert .  \label{2.20}
\end{equation}%
The equality holds in (\ref{2.20}) if and only if%
\begin{equation}
\int_{a}^{b}f\left( t\right) dt=2\left( \frac{\sqrt{m_{1}M_{1}}}{M_{1}+m_{1}}%
+i\frac{\sqrt{m_{2}M_{2}}}{M_{2}+m_{2}}\right) \left( \int_{a}^{b}\left\Vert
f\left( t\right) \right\Vert dt\right) e.  \label{2.21}
\end{equation}
\end{corollary}

\begin{proof}
Firstly, remark that, for $x,z,Z\in H,$ the following statements are
equivalent.

\begin{enumerate}
\item[(i)] $\func{Re}\left\langle Z-x,x-z\right\rangle \geq 0$

and

\item[(ii)] $\left\Vert x-\frac{Z+z}{2}\right\Vert \leq \frac{1}{2}%
\left\Vert Z-z\right\Vert .$
\end{enumerate}

Using this fact, we may simply realise that (\ref{2.18}) and (\ref{2.19})
are equivalent.

Now, from the first inequality in (\ref{2.18}), we get%
\begin{equation*}
\left\Vert f\left( t\right) \right\Vert ^{2}+m_{1}M_{1}\leq \left(
M_{1}+m_{1}\right) \func{Re}\left\langle f\left( t\right) ,e\right\rangle
\end{equation*}%
implying%
\begin{equation}
\frac{\left\Vert f\left( t\right) \right\Vert ^{2}}{\sqrt{m_{1}M_{1}}}+\sqrt{%
m_{1}M_{1}}\leq \frac{M_{1}+m_{1}}{\sqrt{m_{1}M_{1}}}\func{Re}\left\langle
f\left( t\right) ,e\right\rangle  \label{2.22}
\end{equation}%
for a.e. $t\in \left[ a,b\right] .$

Since, obviously,%
\begin{equation}
2\left\Vert f\left( t\right) \right\Vert \leq \frac{\left\Vert f\left(
t\right) \right\Vert ^{2}}{\sqrt{m_{1}M_{1}}}+\sqrt{m_{1}M_{1}},
\label{2.23}
\end{equation}%
hence, by (\ref{2.22}) and (\ref{2.23})%
\begin{equation}
0\leq \frac{2\sqrt{m_{1}M_{1}}}{M_{1}+m_{1}}\left\Vert f\left( t\right)
\right\Vert \leq \func{Re}\left\langle f\left( t\right) ,e\right\rangle
\label{2.24}
\end{equation}%
for a.e. $t\in \left[ a,b\right] .$

Using the same argument as in the proof of Corollary \ref{c2.2}, we deduce
the desired inequality. We omit the details.
\end{proof}

\section{The Case of Orthonormal Vectors}

In the early paper \cite{SSD1}, we pointed out the following reverse of the
continuous triangle inequality for real or complex Hilbert spaces $\left(
H;\left\langle \cdot ,\cdot \right\rangle \right) .$

\begin{theorem}
\label{t3.1}Let $\left\{ e_{1},\dots ,e_{n}\right\} $ be a family of
orthonormal vectors in $H$, $k_{i}\geq 0,$ $i\in \left\{ 1,\dots ,n\right\} $
and $f\in L\left( \left[ a,b\right] ;H\right) $ such that 
\begin{equation}
k_{i}\left\Vert f\left( t\right) \right\Vert \leq \func{Re}\left\langle
f\left( t\right) ,e_{i}\right\rangle   \label{3.1}
\end{equation}%
for each $i\in \left\{ 1,\dots ,n\right\} $ and for a.e. $t\in \left[ a,b%
\right] .$ Then%
\begin{equation}
\left( \sum_{i=1}^{n}k_{i}^{2}\right) ^{\frac{1}{2}}\int_{a}^{b}\left\Vert
f\left( t\right) \right\Vert dt\leq \left\Vert \int_{a}^{b}f\left( t\right)
dt\right\Vert ,  \label{3.2}
\end{equation}%
where the case of equality holds if and only if%
\begin{equation}
\int_{a}^{b}f\left( t\right) dt=\left( \int_{a}^{b}\left\Vert f\left(
t\right) \right\Vert dt\right) \sum_{i=1}^{n}k_{i}e_{i}.  \label{3.3}
\end{equation}
\end{theorem}

In what follows, we improve this result for the case of complex Hilbert
spaces. The following result holds.

\begin{theorem}
\label{t3.2}Let $\left\{ e_{1},\dots ,e_{n}\right\} $ be a family of
orthonormal vectors in the complex Hilbert space $\left( H;\left\langle
\cdot ,\cdot \right\rangle \right) $. If $k_{j},h_{j}\geq 0,$ $j\in \left\{
1,\dots ,n\right\} $ and $f\in L\left( \left[ a,b\right] ;H\right) $ are
such that%
\begin{equation}
k_{j}\left\Vert f\left( t\right) \right\Vert \leq \func{Re}\left\langle
f\left( t\right) ,e_{j}\right\rangle ,\quad h_{j}\left\Vert f\left( t\right)
\right\Vert \leq \func{Im}\left\langle f\left( t\right) ,e_{j}\right\rangle 
\label{3.4}
\end{equation}%
for each $j\in \left\{ 1,\dots ,n\right\} $ and a.e. $t\in \left[ a,b\right]
,$ then%
\begin{equation}
\left[ \sum_{j=1}^{n}\left( k_{j}^{2}+h_{j}^{2}\right) \right] ^{\frac{1}{2}%
}\int_{a}^{b}\left\Vert f\left( t\right) \right\Vert dt\leq \left\Vert
\int_{a}^{b}f\left( t\right) dt\right\Vert .  \label{3.5}
\end{equation}%
The case of equality holds in (\ref{3.5}) if and only if%
\begin{equation}
\int_{a}^{b}f\left( t\right) dt=\left( \int_{a}^{b}\left\Vert f\left(
t\right) \right\Vert dt\right) \sum_{j=1}^{n}\left( k_{j}+ih_{j}\right)
e_{j}.  \label{3.6}
\end{equation}
\end{theorem}

\begin{proof}
Before we prove the theorem, let us recall that, if $x\in H$ and $%
e_{1},\dots ,e_{n}$ are orthonormal vectors, then the following identity
holds true:%
\begin{equation}
\left\Vert x-\sum_{j=1}^{n}\left\langle x,e_{j}\right\rangle
e_{j}\right\Vert ^{2}=\left\Vert x\right\Vert ^{2}-\sum_{j=1}^{n}\left\vert
\left\langle x,e_{j}\right\rangle \right\vert ^{2}.  \label{3.7}
\end{equation}%
As a consequence of this identity, we have the \textit{Bessel inequality}%
\begin{equation}
\sum_{j=1}^{n}\left\vert \left\langle x,e_{j}\right\rangle \right\vert
^{2}\leq \left\Vert x\right\Vert ^{2},x\in H.  \label{3.8}
\end{equation}%
in which the case of equality holds if and only if 
\begin{equation}
x=\sum_{j=1}^{n}\left\langle x,e_{j}\right\rangle e_{j}.  \label{3.9}
\end{equation}%
Now, applying Bessel's inequality for $x=\int_{a}^{b}f\left( t\right) dt,$
we have successively%
\begin{align}
\left\Vert \int_{a}^{b}f\left( t\right) dt\right\Vert ^{2}& \geq
\sum_{j=1}^{n}\left\vert \left\langle \int_{a}^{b}f\left( t\right)
dt,e_{j}\right\rangle \right\vert ^{2}=\sum_{j=1}^{n}\left\vert
\int_{a}^{b}\left\langle f\left( t\right) ,e_{j}\right\rangle dt\right\vert
^{2}  \label{3.10} \\
& =\sum_{j=1}^{n}\left\vert \int_{a}^{b}\func{Re}\left\langle f\left(
t\right) ,e_{j}\right\rangle dt+i\left( \int_{a}^{b}\func{Im}\left\langle
f\left( t\right) ,e_{j}\right\rangle dt\right) \right\vert ^{2}  \notag \\
& =\sum_{j=1}^{n}\left[ \left( \int_{a}^{b}\func{Re}\left\langle f\left(
t\right) ,e_{j}\right\rangle dt\right) ^{2}+\left( \int_{a}^{b}\func{Im}%
\left\langle f\left( t\right) ,e_{j}\right\rangle dt\right) ^{2}\right] . 
\notag
\end{align}%
Integrating (\ref{3.4}) on $\left[ a,b\right] ,$ we get%
\begin{equation}
\int_{a}^{b}\func{Re}\left\langle f\left( t\right) ,e_{j}\right\rangle
dt\geq k_{j}\int_{a}^{b}\left\Vert f\left( t\right) \right\Vert dt
\label{3.11}
\end{equation}%
and%
\begin{equation}
\int_{a}^{b}\func{Im}\left\langle f\left( t\right) ,e_{j}\right\rangle
dt\geq h_{j}\int_{a}^{b}\left\Vert f\left( t\right) \right\Vert dt
\label{3.12}
\end{equation}%
for each $j\in \left\{ 1,\dots ,n\right\} .$

Squaring and adding the above two inequalities (\ref{3.11}) and (\ref{3.12}%
), we deduce%
\begin{multline*}
\sum_{j=1}^{n}\left[ \left( \int_{a}^{b}\func{Re}\left\langle f\left(
t\right) ,e_{j}\right\rangle dt\right) ^{2}+\left( \int_{a}^{b}\func{Im}%
\left\langle f\left( t\right) ,e_{j}\right\rangle dt\right) ^{2}\right] \\
\geq \sum_{j=1}^{n}\left( k_{j}^{2}+h_{j}^{2}\right) \left(
\int_{a}^{b}\left\Vert f\left( t\right) \right\Vert dt\right) ^{2},
\end{multline*}%
which combined with (\ref{3.10}) will produce the desired inequality (\ref%
{3.5}).

Now, if (\ref{3.6}) holds true, then%
\begin{align*}
\left\Vert \int_{a}^{b}f\left( t\right) dt\right\Vert & =\left(
\int_{a}^{b}\left\Vert f\left( t\right) \right\Vert dt\right) \left\Vert
\sum_{j=1}^{n}\left( k_{j}+ih_{j}\right) e_{j}\right\Vert \\
& =\left( \int_{a}^{b}\left\Vert f\left( t\right) \right\Vert dt\right)
\left( \left\Vert \sum_{j=1}^{n}\left( k_{j}+ih_{j}\right) e_{j}\right\Vert
^{2}\right) ^{\frac{1}{2}} \\
& =\left( \int_{a}^{b}\left\Vert f\left( t\right) \right\Vert dt\right) 
\left[ \sum_{j=1}^{n}\left( k_{j}^{2}+h_{j}^{2}\right) \right] ^{\frac{1}{2}%
},
\end{align*}%
and the case of equality holds in (\ref{3.5}).

Conversely, if the equality holds in (\ref{3.5}), then it must hold in all
the inequalities used to prove (\ref{3.5}) and therefore we must have%
\begin{equation}
\left\Vert \int_{a}^{b}f\left( t\right) dt\right\Vert
^{2}=\sum_{j=1}^{n}\left\vert \left\langle \int_{a}^{b}f\left( t\right)
dt,e_{j}\right\rangle \right\vert ^{2}  \label{3.13}
\end{equation}%
and%
\begin{equation}
k_{j}\left\Vert f\left( t\right) \right\Vert =\func{Re}\left\langle f\left(
t\right) ,e_{j}\right\rangle \text{ \ and \ }h_{j}\left\Vert f\left(
t\right) \right\Vert =\func{Re}\left\langle f\left( t\right)
,e_{j}\right\rangle  \label{3.14}
\end{equation}%
for each $j\in \left\{ 1,\dots ,n\right\} $ and a.e. $t\in \left[ a,b\right]
.$

From (\ref{3.13}), on using the identity (\ref{3.9}), we deduce that%
\begin{equation}
\int_{a}^{b}f\left( t\right) dt=\sum_{j=1}^{n}\left\langle
\int_{a}^{b}f\left( t\right) dt,e_{j}\right\rangle e_{j}.  \label{3.15}
\end{equation}%
Now, multiplying the second equality in (\ref{3.14}) with the imaginary unit 
$i,$ integrating both inequalities on $\left[ a,b\right] $ and summing them
up, we get%
\begin{equation}
\left( k_{j}+ih_{j}\right) \int_{a}^{b}\left\Vert f\left( t\right)
\right\Vert dt=\left\langle \int_{a}^{b}f\left( t\right)
dt,e_{j}\right\rangle  \label{3.16}
\end{equation}%
for each $j\in \left\{ 1,\dots ,n\right\} .$

Finally, utilising (\ref{3.15}) and (\ref{3.16}), we deduce (\ref{3.6}) and
the theorem is proved.
\end{proof}

The following corollaries are of interest.

\begin{corollary}
\label{c3.2}Let $e_{1},\dots ,e_{m}$ be orthonormal vectors in the complex
Hilbert space $\left( H;\left\langle \cdot ,\cdot \right\rangle \right) $
and $\rho _{k},\eta _{k}\in \left( 0,1\right) ,$ $k\in \left\{ 1,\dots
,n\right\} .$ If $f\in L\left( \left[ a,b\right] ;H\right) $ is such that%
\begin{equation*}
\left\Vert f\left( t\right) -e_{k}\right\Vert \leq \rho _{k},\qquad
\left\Vert f\left( t\right) -ie_{k}\right\Vert \leq \eta _{k}
\end{equation*}%
for each $k\in \left\{ 1,\dots ,n\right\} $ and for a.e. $t\in \left[ a,b%
\right] ,$ then we have the inequality%
\begin{equation}
\left[ \sum_{k=1}^{n}\left( 2-\rho _{k}^{2}-\eta _{k}^{2}\right) \right] ^{%
\frac{1}{2}}\int_{a}^{b}\left\Vert f\left( t\right) \right\Vert dt\leq
\left\Vert \int_{a}^{b}f\left( t\right) dt\right\Vert .  \label{3.17}
\end{equation}%
The case of equality holds in (\ref{3.17}) if and only if%
\begin{equation}
\int_{a}^{b}f\left( t\right) dt=\left( \int_{a}^{b}\left\Vert f\left(
t\right) \right\Vert dt\right) \sum_{k=1}^{n}\left( \sqrt{1-\rho _{k}^{2}}+i%
\sqrt{1-\eta _{k}^{2}}\right) e_{k}.  \label{3.18}
\end{equation}
\end{corollary}

The proof follows by Theorem \ref{t3.2} and is similar to the one from
Corollary \ref{c2.2}. We omit the details.

\begin{corollary}
\label{c3.3}Let $e_{1},\dots ,e_{m}$ be as in Corollary \ref{c3.2} and $%
M_{k}\geq m_{k}>0,$ $N_{k}\geq n_{k}>0,$ $k\in \left\{ 1,\dots ,n\right\} .$
If $f\in L\left( \left[ a,b\right] ;H\right) $ is such that either%
\begin{equation*}
\func{Re}\left\langle M_{k}e_{k}-f\left( t\right) ,f\left( t\right)
-m_{k}e_{k}\right\rangle \geq 0,\ \ \func{Re}\left\langle
N_{k}ie_{k}-f\left( t\right) ,f\left( t\right) -n_{k}ie_{k}\right\rangle
\geq 0
\end{equation*}%
or, equivalently,%
\begin{align*}
\left\Vert f\left( t\right) -\frac{M_{k}+m_{k}}{2}e_{k}\right\Vert & \leq 
\frac{1}{2}\left( M_{k}-m_{k}\right) ,\  \\
\left\Vert f\left( t\right) -\frac{N_{k}+n_{k}}{2}ie_{k}\right\Vert & \leq 
\frac{1}{2}\left( N_{k}-n_{k}\right)
\end{align*}%
for each $k\in \left\{ 1,\dots ,n\right\} $ and a.e. $t\in \left[ a,b\right]
,$ then we have the inequality%
\begin{equation}
2\left\{ \sum_{k=1}^{m}\left[ \frac{m_{k}M_{k}}{\left( M_{k}+m_{k}\right)
^{2}}+i\frac{n_{k}N_{k}}{\left( N_{k}+n_{k}\right) ^{2}}\right] \right\} ^{%
\frac{1}{2}}\int_{a}^{b}\left\Vert f\left( t\right) \right\Vert dt\leq
\left\Vert \int_{a}^{b}f\left( t\right) dt\right\Vert .  \label{3.19}
\end{equation}%
The case of equality holds in (\ref{3.19}) if and only if%
\begin{equation}
\int_{a}^{b}f\left( t\right) dt=2\left( \int_{a}^{b}\left\Vert f\left(
t\right) \right\Vert dt\right) \sum_{k=1}^{n}\left( \frac{\sqrt{m_{k}M_{k}}}{%
M_{k}+m_{k}}+i\frac{\sqrt{n_{k}N_{k}}}{N_{k}+n_{k}}\right) e_{k}.
\label{3.20}
\end{equation}
\end{corollary}

The proof employs Theorem \ref{t3.2} and is similar to the one in Corollary %
\ref{c2.3}. We omit the details.

\section{Applications for Complex-Valued Functions}

The following reverse of the generalised triangle inequality for
complex-valued functions that improves Karamata's result (\ref{1.2}) holds.

\begin{proposition}
\label{p4.1}Let $f\in L\left( \left[ a,b\right] ;\mathbb{C}\right) $ with
the property that%
\begin{equation}
0\leq \varphi _{1}\leq \arg f\left( t\right) \leq \varphi _{2}<\frac{\pi }{2}
\label{4.1}
\end{equation}%
for a.e. $t\in \left[ a,b\right] .$ Then we have the inequality%
\begin{equation}
\sqrt{\sin ^{2}\varphi _{1}+\cos ^{2}\varphi _{2}}\int_{a}^{b}\left\vert
f\left( t\right) \right\vert dt\leq \left\vert \int_{a}^{b}f\left( t\right)
dt\right\vert .  \label{4.2}
\end{equation}%
The equality holds in (\ref{4.2}) if and only if%
\begin{equation}
\int_{a}^{b}f\left( t\right) dt=\left( \cos \varphi _{2}+i\sin \varphi
_{1}\right) \int_{a}^{b}\left\vert f\left( t\right) \right\vert dt.
\label{4.3}
\end{equation}
\end{proposition}

\begin{proof}
Let $f\left( t\right) =\func{Re}f\left( t\right) +i\func{Im}f\left( t\right)
.$ We may assume that $\func{Re}f\left( t\right) \geq 0,$ $\func{Im}f\left(
t\right) >0,$ for a.e. $t\in \left[ a,b\right] ,$ since, by (\ref{4.1}), $%
\frac{\func{Im}f\left( t\right) }{\func{Re}f\left( t\right) }=\tan \left[
\arg f\left( t\right) \right] \in \left[ 0,\frac{\pi }{2}\right) $, for a.e. 
$t\in \left[ a,b\right] .$ By (\ref{4.1}), we obviously have%
\begin{equation*}
0\leq \tan ^{2}\varphi _{1}\leq \left[ \frac{\func{Im}f\left( t\right) }{%
\func{Re}f\left( t\right) }\right] ^{2}\leq \tan ^{2}\varphi _{2},\qquad 
\text{for a.e. }t\in \left[ a,b\right] ,
\end{equation*}%
from where we get%
\begin{equation*}
\frac{\left[ \func{Im}f\left( t\right) \right] ^{2}+\left[ \func{Re}f\left(
t\right) \right] ^{2}}{\left[ \func{Re}f\left( t\right) \right] ^{2}}\leq 
\frac{1}{\cos {}^{2}\varphi _{2}},
\end{equation*}%
for a.e. $t\in \left[ a,b\right] ,$ and%
\begin{equation*}
\frac{\left[ \func{Im}f\left( t\right) \right] ^{2}+\left[ \func{Re}f\left(
t\right) \right] ^{2}}{\left[ \func{Im}f\left( t\right) \right] ^{2}}\leq 
\frac{1+\tan ^{2}\varphi _{1}}{\tan ^{2}\varphi _{1}}=\frac{1}{\sin \varphi
_{1}},
\end{equation*}%
for a.e. $t\in \left[ a,b\right] ,$ giving the simpler inequalities%
\begin{equation*}
\left\vert f\left( t\right) \right\vert \cos \varphi _{2}\leq \func{Re}%
\left( f\left( t\right) \right) ,\quad \left\vert f\left( t\right)
\right\vert \sin \varphi _{1}\leq \func{Im}\left( f\left( t\right) \right)
\end{equation*}%
for a.e. $t\in \left[ a,b\right] .$

Now, applying Theorem \ref{t2.1} for the complex Hilbert space $\mathbb{C}$
endowed with the inner product $\left\langle z,w\right\rangle =z\cdot \bar{w}
$ for $k_{1}=\cos \varphi _{2},$ $k_{2}=\sin \varphi _{1}$ and $e=1,$ we
deduce the desired inequality (\ref{4.2}). The case of equality is also
obvious and we omit the details.
\end{proof}

Another result that has an obvious geometrical interpretation is the
following one.

\begin{proposition}
\label{p4.2}Let $e\in \mathbb{C}$ with $\left\vert e\right\vert =1$ and $%
\rho _{1},\rho _{2}\in \left( 0,1\right) .$ If $f\left( t\right) \in L\left( %
\left[ a,b\right] ;\mathbb{C}\right) $ such that%
\begin{equation}
\left\vert f\left( t\right) -e\right\vert \leq \rho _{1},\quad \left\vert
f\left( t\right) -ie\right\vert \leq \rho _{2}\qquad \text{for a.e. }t\in %
\left[ a,b\right] ,  \label{4.5}
\end{equation}%
then we have the inequality%
\begin{equation}
\sqrt{2-\rho _{1}^{2}-\rho _{2}^{2}}\int_{a}^{b}\left\vert f\left( t\right)
\right\vert dt\leq \left\vert \int_{a}^{b}f\left( t\right) dt\right\vert ,
\label{4.6}
\end{equation}%
with equality if and only if%
\begin{equation}
\int_{a}^{b}f\left( t\right) dt=\left( \sqrt{1-\rho _{1}^{2}}+i\sqrt{1-\rho
_{2}^{2}}\right) \int_{a}^{b}\left\vert f\left( t\right) \right\vert dt\cdot
e.  \label{4.7}
\end{equation}
\end{proposition}

The proof is obvious by Corollary \ref{c2.2} applied for $H=\mathbb{C}$ and
we omit the details.

\begin{remark}
If we choose $e=1,$ and for $\rho _{1},\rho _{2}\in \left( 0,1\right) $ we
define 
\begin{equation*}
\bar{D}\left( 1,\rho _{1}\right) :=\left\{ z\in \mathbb{C}|\left\vert
z-1\right\vert \leq \rho _{1}\right\} ,\ \ \bar{D}\left( i,\rho _{2}\right)
:=\left\{ z\in \mathbb{C}|\left\vert z-i\right\vert \leq \rho _{2}\right\} ,
\end{equation*}%
then obviously the intersection domain%
\begin{equation*}
S_{\rho _{1},\rho _{2}}:=\bar{D}\left( 1,\rho _{1}\right) \cap \bar{D}\left(
i,\rho _{2}\right)
\end{equation*}%
is nonempty if and only if $\rho _{1}+\rho _{2}>\sqrt{2}.$

If $f\left( t\right) \in S_{\rho _{1},\rho _{2}}$ for a.e. $t\in \left[ a,b%
\right] ,$ then (\ref{4.6}) holds true. The equality holds in (\ref{4.6}) if
and only if%
\begin{equation*}
\int_{a}^{b}f\left( t\right) dt=\left( \sqrt{1-\rho _{1}^{2}}+i\sqrt{1-\rho
_{2}^{2}}\right) \int_{a}^{b}\left\vert f\left( t\right) \right\vert dt.
\end{equation*}
\end{remark}

\end{document}